\documentclass[12pt]{article}
\usepackage{amsfonts}
\usepackage{a4wide}
\usepackage{amsmath,amscd,amsthm,a4,amssymb}

\begin{document}

\begin{center}
{\large \textbf{ON THE REPRESENTATION BY LINEAR SUPERPOSITIONS
}}\footnote{ This research has been supported by INTAS under Grant
YSF-06-1000015-6283}

\

\textbf{Vugar E. Ismailov} \

{Mathematics and Mechanics Institute}

Azerbaijan National Academy of Sciences

Az-1141, Baku, Azerbaijan

{e-mail:} {vugaris@mail.ru}
\end{center}

\textbf{Abstract.} In a number of papers, Y. Sternfeld investigated the
problems of representation of continuous and bounded functions by linear
superpositions. In particular, he proved that if such representation holds
for continuous functions, then it holds for bounded functions. We consider
the same problem without involving any topology and establish a rather
practical necessary and sufficient condition for representability of an
arbitrary function by linear superpositions. In particular, we show that if
some representation by linear superpositions holds for continuous functions,
then it holds for all functions. This will lead us to the analogue of the
well-known Kolmogorov superposition theorem for multivariate functions on
the $d$-dimensional unit cube.

\bigskip



\bigskip

\begin{center}
{\large \textbf{1. Introduction}}
\end{center}

Let $X$ be any set and $h_{i}:X\rightarrow {{\mathbb{R}}},~i=1,...,r,$ be
arbitrarily fixed functions. Consider the set
\begin{equation*}
\mathcal{B}(X)=\mathcal{B}(h_{1},...,h_{r};X)=\left\{
\sum\limits_{i=1}^{r}g_{i}(h_{i}(x)),~x\in X,~g_{i}:\mathbb{R}\rightarrow
\mathbb{R},~i=1,...,r\right\} \eqno(1.1)
\end{equation*}
Members of this set will be called linear superpositions (see
[35]). We are going to answer the question: what conditions on $X$
guarantee that each function on $X$ will be in the set
$\mathcal{B}(X)$? The simplest case $ X\subset
\mathbb{R}^{d},~r=d$ and $h_{i}$ are the coordinate functions has
been solved in [16]. See also [5,15] for the case $r=2.$

By $\mathcal{B}_{c}(X)$ and $\mathcal{B}_{b}(X)$ denote the right
hand side of (1.1) with continuous and bounded
$g_{i}:\mathbb{R}\rightarrow \mathbb{R} ,~i=1,...,r,$
respectively. Our starting point is the well-known superposition
theorem of Kolmogorov [17]. It states that for the unit cube $
\mathbb{I}^{d},~\mathbb{I}=[0,1],~d\geq 2,$ there exists $2d+1$
functions $ \{s_{q}\}_{q=1}^{2d+1}\subset C(\mathbb{I}^{d})$ of
the form

\begin{equation*}
s_{q}(x_{1},...,x_{d})=\sum_{p=1}^{d}\varphi _{pq}(x_{p}),~\varphi _{pq}\in
C(\mathbb{I}),~p=1,...,d,~q=1,...,2d+1\eqno(1.2)
\end{equation*}
such that each function $f\in C(\mathbb{I}^{d})$ admits the representation

\begin{equation*}
f(x)=\sum_{q=1}^{2d+1}g_{q}(s_{q}(x)),~x=(x_{1},...,x_{d})\in
\mathbb{I} ^{d},~g_{q}\in C({{\mathbb{R)}}}.\eqno(1.3)
\end{equation*}

In our notation, (1.3) means that
$\mathcal{B}_{c}(s_{1},...,s_{2d+1};
\mathbb{I}^{d})=C(\mathbb{I}^{d}).$ This surprising and deep
result, which solved (negatively) Hilbert's 13-th problem, was
improved and generalized in several directions. It was first
observed by Lorentz [21] that the functions $g_{q}$ can be
replaced by a single continuous function $g.$ Sprecher [29] showed
that the theorem can be proven with constant multiples of a single
function $\varphi $ and translations. Specifically, $\varphi
_{pq}$ in (1.2) can be chosen as $\lambda ^{p}\varphi
(x_{p}+\varepsilon q),$ where $ \varepsilon $ and $\lambda $ are
some positive constants. Fridman [9] succeeded in showing that the
functions $\varphi _{pq}$ can be constructed to belong to the
class $Lip(1).$ Vitushkin and Henkin [35] showed that $ \varphi
_{pq}$ cannot be taken to be continuously differentiable.

Ostrand [25] extended the Kolmogorov theorem to general compact
metric spaces. In particular, he proved that for each compact
$d$-dimensional metric space $X$ there exist continuous real
functions $\{\alpha _{i}\}_{i=1}^{2d+1}\subset C(X)$ such that
$\mathcal{B}_{c}(\alpha _{1},...,\alpha _{2d+1};X)=C(X).$
Sternfeld [32] showed that the number $ 2d+1 $ cannot be reduced
for any $d$-dimensional space $X.$ Thus the number of terms in the
Kolmogorov superposition theorem is the best possible.

Some papers of Sternfeld have been devoted to the representation
of continuous and bounded functions by linear superpositions. Let
$C(X)$ and $ B(X)$ denote the space of continuous and bounded
functions on some set $X$ respectively (in the first case, $X$ is
supposed to be a compact metric space). Let $F=\{h\}$ be a family
of functions on $X.$ $F$ is called a uniformly separating family
(\textit{u.s.f.}) if there exists a number $ 0<\lambda \leq 1$
such that for each pair $\{x_{j}\}_{j=1}^{m}$, $
\{z_{j}\}_{j=1}^{m}$ of disjoint finite sequences in $X$, there
exists some $ h\in F$ so that if from the two sequences
$\{h(x_{j})\}_{j=1}^{m}$and $ \{h(z_{j})\}_{j=1}^{m}$ in $h(X)$ we
remove a maximal number of pairs of points $h(x_{j_{1}})$ and
$h(z_{j_{2}})$ with $h(x_{j_{1}})=h(z_{j_{2}}),$ there remains at
least $\lambda m$ points in each sequence (or , equivalently, at
most $(1-\lambda )m$ pairs can be removed). Sternfeld [31] proved
that for a finite family $F=\{h_{1},...,h_{r}\}$ of functions on
$X$, being a \textit{u.s.f.} is equivalent to the equality
$\mathcal{B} _{b}(h_{1},...,h_{r};X)=B(X),$ and that in the case
where $X$ is a compact metric space and the elements of $F$ are
continuous functions on $X$, the equality
$\mathcal{B}_{c}(h_{1},...,h_{r};X)=C(X)$ implies that $F$ is a
\textit{u.s.f.} Thus, in particular, Sternfeld obtained that the
formula (1.3) is valid for all bounded functions, where $g_{q}$
are bounded functions depending on $f$ (see also [15, p.21]).

Let $X$ be a compact metric space. The family $F=\{h\}\subset
C(X)$ is said to be a measure separating family (\textit{m.s.f.})
if there exists a number $0<\lambda \leq 1$ such that for any
measure $\mu $ in $\ C(X)^{\ast },$ the inequality $\left\Vert \mu
\circ h^{-1}\right\Vert \geq \lambda \left\Vert \mu \right\Vert $
holds for some $h\in F.$ Sternfeld [33] proved that $
\mathcal{B}_{c}(h_{1},...,h_{r};X)=C(X)$ if and only if the family
$ \{h_{1},...,h_{r}\}$ is a \textit{m.s.f.} In [31], it has been
shown that if $r=2,$ then the properties \textit{u.s.f.} and
\textit{m.s.f.} are equivalent. Therefore, the equality
$\mathcal{B}_{b}(h_{1},h_{2};X)=B(X)$ is equivalent to
$\mathcal{B}_{c}(h_{1},h_{2};X)=C(X).$ But for $r\,>2$, these two
properties are no longer equivalent. That is, $\mathcal{B}
_{b}(h_{1},...,h_{r};X)=B(X)$ does not always imply $\mathcal{B}
_{c}(h_{1},...,h_{r};X)=C(X)$ (see [33]).

Our purpose is to consider the above mentioned problem of
representation by linear superpositions without involving any
topology (that of continuity or boundedness). More precisely, we
want to characterize those sets $X$ for which
$\mathcal{B}(h_{1},...,h_{r};X)=T(X),$ where $T(X)$ is the space
of all functions on $X.$ This will be done in terms of
\textit{closed paths}, the explicit and practically convenient
objects. We show that nonexistence of closed paths in $X$ is
equivalent to the equality $\mathcal{B}(X)=T(X).$ In particular,
we will obtain that $\mathcal{B}_{c}(X)=C(X)$ implies $
\mathcal{B}(X)=T(X).$ Therefore, the formula (1.3) is valid for
all multivariate functions over the unite cube $\mathbb{I}^{d},$
where $g_{q}$ are univariate functions depending on $f.$ We will
also present an example due to Khavinson [15] showing that even in
the case $r=2,$ the equality $ \mathcal{B}(h_{1},h_{2};X)=T(X)$
does not imply $\mathcal{B} _{c}(h_{1},h_{2};X)=C(X).$ At the end
we will make some observations around the problems of
representation and interpolation by \textit{ridge} functions,
which are widely used in multivariate approximation theory.

\bigskip

\begin{center}
{\large \textbf{2. Closed paths}}
\end{center}

In the sequel, by $\delta _{A}$ we will denote the characteristic function
of a set $\ A\subset \mathbb{R}.$ That is,

\begin{equation*}
\delta _{A}(y)=\left\{
\begin{array}{c}
1,~if~y\in A \\ 0,~if~y\notin A.
\end{array}
\right.
\end{equation*}

The following definition is based on the ideas set forth in the works [1]
and [16].

\bigskip

\textbf{Definition 2.1.} \textit{Given a set $X$ and functions $
h_{i}:X\rightarrow \mathbb{R},~i=1,...,r$. A set of points $
\{x_{1},...,x_{n}\}\subset X$ is called to be a closed path with
respect to the functions $h_{1},...,h_{r}$ (or, concisely, a
closed path if there is no confusion), if there exists a vector
$\lambda =(\lambda _{1},...,\lambda _{n})$ with the nonzero real
coordinates $\lambda _{i},~i=1,...,n,$ such that }

\begin{equation*}
\sum_{j=1}^{n}\lambda _{j}\delta _{h_{i}(x_{j})}=0,~i=1,...,r.\eqno(2.1)
\end{equation*}

\bigskip

Let for $i=1,...,r,$ the set $\{h_{i}(x_{j}),~j=1,...,n\}$ have
$k_{i}$ different values. Then it is not difficult to see that Eq.
(2.1) stands for a system of $\sum_{i=1}^{r}k_{i}$ homogeneous
linear equations in unknowns $ \lambda _{1},...,\lambda _{n}.$ If
this system has any solution with the nonzero components, then the
given set $\{x_{1},...,x_{n}\}$ is a closed path. In the last
case, the system has also a solution $m=(m_{1},...,m_{n})$ with
the nonzero integer components $m_{i},~i=1,...,n.$ Thus, in
Definition 2.1, the vector $\lambda =(\lambda _{1},...,\lambda
_{n})$ can be replaced by a vector $m=(m_{1},...,m_{n})$ with
$m_{i}\in \mathbb{Z}\backslash \{0\}.$

For example, the set
$l=\{(0,0,0),~(0,0,1),~(0,1,0),~(1,0,0),~(1,1,1)\}$ is a closed
path in $\mathbb{R}^{3}$ with respect to the functions $
h_{i}(z_{1},z_{2},z_{3})=z_{i},~i=1,2,3.$ The vector $\lambda $ in
Definition 2.1 can be taken as $(-2,1,1,1,-1).$

In the case $r=2,$ the picture of closed path becomes more clear.
Let, for example, $h_{1}$ and $h_{2}$ be the coordinate functions
on $\mathbb{R}^{2}.$ In this case, a closed path is the union of
some sets $A_{k}$ with the property: each $A_{k}$ consists of
vertices of a closed broken line with the sides parallel to the
coordinate axis. These objects (sets $A_{k}$) have been exploited
in practically all works devoted to the approximation of bivariate
functions by univariate functions, although under the different
names (see, for example, [15, chapter 2]). If $X$ and the
functions $h_{1}$ and $h_{2}$ are arbitrary, the sets $A_{k}$ can
be described as a trace of some point traveling alternatively in
the level sets of $h_{1}$ and $h_{2},$ and then returning to its
primary position. It should be remarked that in the case $r>2,$
closed paths do not admit such a simple geometric description. We
refer the reader to Braess and Pinkus [1] for the description of
closed paths when $r=3$ and $h_{i}(\mathbf{x})=\mathbf{a}^{i}
\mathbf{\cdot x},$ $\mathbf{x\in
}\mathbb{R}^{2},~\mathbf{a}^{i}\in \mathbb{R }^{2}\backslash
\{\mathbf{0}\},~i=1,2,3.$

Let $T(X)$ denote the set of all functions on $X.$ With each pair
$ \left\langle p,\lambda \right\rangle ,$ where
$p=\{x_{1},...,x_{n}\}$ is a closed path in $X$ and $\lambda
=(\lambda _{1},...,\lambda _{n})$ is a vector known from
Definition 2.1, we associate the functional

\begin{equation*}
G_{p,\lambda }:T(X)\rightarrow \mathbb{R},~~G_{p,\lambda
}(f)=\sum_{j=1}^{n}\lambda _{j}f(x_{j}).
\end{equation*}
In the following, such pairs $\left\langle p,\lambda \right\rangle
$ will be called \textit{closed path-vector pairs} of $X.$ It is
clear that the functional $G_{p,\lambda }$ is linear. Besides,
$G_{p,\lambda }(g)=0$ for all functions $g\in
\mathcal{B}(h_{1},...,h_{r};X).$ Indeed, assume that (2.1) holds.
Given $i\leq r$, let $z=h_{i}(x_{j})$ for some $j$. Hence, $
\sum_{j~(h_{i}(x_{j})=z)}\lambda _{j}=0$ and
$\sum_{j~(h_{i}(x_{j})=z)} \lambda _{j}g_{i}(h_{i}(x_{j}))=0$. A
summation yields $G_{p,\lambda }(g_{i}\circ h_{i})=0$. Since
$G_{p,\lambda }$ is linear, we obtain that $ G_{p,\lambda
}(\sum_{i=1}^{r}g_{i}\circ h_{i})=0$.

The following lemma characterizes the set $\mathcal{B}(h_{1},...,h_{r};X)$
under some restrictions and will be used in the proof of Theorem 2.5 below.

\bigskip

\textbf{Lemma 2.2.} \textit{Let $X$ have closed paths
$h_{i}(X)\cap h_{j}(X)=\varnothing ,$ for all $i,j\in
\{1,...,r\},~i\neq j.$ Then a function $f:X\rightarrow \mathbb{R}$
belongs to the set $\mathcal{B} (h_{1},...,h_{r};X)$ if and only
if $G_{p,\lambda }(f)=0$ for any closed path-vector pair
$\left\langle p,\lambda \right\rangle $ of $X.$}

\begin{proof} The necessity is obvious, since the functional $G_{p,\lambda
}$ annihilates all members of \\ $\mathcal{B}(h_{1},...,h_{r};X)$.
Let us prove the sufficiency. Introduce the notation
\begin{eqnarray*}
Y_{i} &=&h_{i}(X),~i=1,...,r; \\
\Omega &=&Y_{1}\cup ...\cup Y_{r}.
\end{eqnarray*}

Consider the following subsets of $\Omega $:
\begin{equation*}
\mathcal{L}=\{Y=\{y_{1},...,y_{r}\}:\text{if there exists }x\in X\text{ such
that }h_{i}(x)=y_{i},~i=1,...,r\}\eqno(2.2)
\end{equation*}

Note that $\mathcal{L}$ is a set of some certain
subsets of $\Omega .$ Each element of $\mathcal{L}$ is a set $Y=\{y_{1},...,y_{r}\}\subset \Omega $ with the property that there exists $x\in X$ such that $h_{i}(x)=y_{i},~i=1,...,r.$ In what follows, all the points $x$ associated with $Y$ by (2.2)
will be called $(\ast )$-points of $Y.$ It is clear that the
number of such points depends on $Y$ as well as on the functions
$h_{1},...,h_{r}$, and may be greater than 1. But note that if any
two points $x_{1}$ and $x_{2}$ are $ (\ast )$-points of $Y$, then
necessarily the set $\{x_{1}$, $x_{2}\}$ forms a closed path with
the associated vector $\lambda =(1;-1).$ In this case, by the
condition of the sufficiency, $f(x_{1})=f(x_{2}).$ Let now
$Y^{\ast }$ be the set of all $(\ast )$-points of $Y.$ Since we
have already known that $ f(Y^{\ast })$ is a single number, we can
define the function

\begin{equation*}
t:\mathcal{L}\rightarrow \mathbb{R},~t(Y)=f(Y^{\ast }).
\end{equation*}
Or, equivalently, $t(Y)=f(x),$ where $x$ is an arbitrary $(\ast
)$-point of $ Y$.

Consider now a class $\mathcal{S}$ of functions of the form $
\sum_{j=1}^{k}r_{j}\delta _{D_{j}},$ where $k$ is a positive
integer, $r_{j}$ are real numbers and $D_{j}$ are elements of
$\mathcal{L},~j=1,...,k.$ We fix neither the numbers $\ k,~r_{j},$
nor the sets $D_{j}.$ Clearly, $ \mathcal{S\ }$is a linear space.
Over $\mathcal{S}$, we define the functional

\begin{equation*}
F:\mathcal{S}\rightarrow \mathbb{R},~F\left( \sum_{j=1}^{k}r_{j}\delta
_{D_{j}}\right) =\sum_{j=1}^{k}r_{j}t(D_{j}).
\end{equation*}

First of all, we must show that this functional is well defined. That is,
the equality

\begin{equation*}
\sum_{j=1}^{k_{1}}r_{j}^{\prime }\delta _{D_{j}^{\prime
}}=\sum_{j=1}^{k_{2}}r_{j}^{\prime \prime }\delta _{D_{j}^{\prime \prime }}
\end{equation*}
always implies the equality

\begin{equation*}
\sum_{j=1}^{k_{1}}r_{j}^{\prime }t(D_{j}^{\prime
})=\sum_{j=1}^{k_{2}}r_{j}^{\prime \prime }t(D_{j}^{\prime \prime }).
\end{equation*}
In fact, this is equivalent to the implication

\begin{equation*}
\sum_{j=1}^{k}r_{j}\delta _{D_{j}}=0\Longrightarrow
\sum_{j=1}^{k}r_{j}t(D_{j})=0,~\text{for all }k\in
\mathbb{N}\text{, } r_{j}\in \mathbb{R}\text{, }D_{j}\subset
\mathcal{L}\text{.}\eqno(2.3)
\end{equation*}

Suppose that the left-hand side of the implication (2.3) be
satisfied. Each set $D_{j}$ consists of $r$ real numbers
$y_{1}^{j},...,y_{r}^{j}$, $ j=1,...,k.$ By the hypothesis of the
lemma, all these numbers are different. Therefore,

\begin{equation*}
\delta _{D_{j}}=\sum_{i=1}^{r}\delta _{y_{i}^{j}},~j=1,...,k.\eqno(2.4)
\end{equation*}
Eq. (2.4) together with the left-hand side of (2.3) gives

\begin{equation*}
\sum_{i=1}^{r}\sum_{j=1}^{k}r_{j}\delta _{y_{i}^{j}}=0.\eqno(2.5)
\end{equation*}
Since the sets $\{y_{i}^{1},y_{i}^{2},...,y_{i}^{k}\}$, $i=1,...,r,$ are
pairwise disjoint, we obtain from (2.5) that

\begin{equation*}
\sum_{j=1}^{k}r_{j}\delta _{y_{i}^{j}}=0,\text{ }i=1,...,r.\eqno(2.6)
\end{equation*}

Let now $x_{1},...,x_{k}$ be some $(\ast )$-points of the sets $
D_{1},...,D_{k}$ respectively. Since by (2.2),
$y_{i}^{j}=h_{i}(x_{j})$, for $i=1,...,r$ and $j=1,...,k,$ it
follows from (2.6) that the set $ \{x_{1},...,x_{k}\}$ is a closed
path. Then by the condition of the sufficiency,
$\sum_{j=1}^{k}r_{j}f(x_{j})=0.$ Hence $
\sum_{j=1}^{k}r_{j}t(D_{j})=0.$ We have proved the implication
(2.3) and hence the functional $F$ is well defined.

Note that the functional $F$ is linear (this can be easily seen from its
definition). Let $F^{\prime }$ be a linear extension of $F$ to the following
space larger than $\mathcal{S}$:

\begin{equation*}
\mathcal{S}^{\prime }=\left\{ \sum_{j=1}^{k}r_{j}\delta _{\omega
_{j}}\right\} ,
\end{equation*}
where $k\in \mathbb{N}$, $r_{j}\in \mathbb{R}$, $\omega
_{j}\subset \Omega .$ As in the above, we do not fix the
parameters $k$, $r_{j}$ and $\omega _{j}.$ Define the functions
\begin{equation*}
g_{i}:Y_{i}\rightarrow \mathbb{R},\text{
}g_{i}(y_{i})\overset{def}{=} F^{\prime }(\delta _{y_{i}}),\text{
}i=1,...,r.
\end{equation*}
Let $x$ be an arbitrary point in $X.$ Obviously, $x$ is a $(\ast )$-point of
some set $Y=\{y_{1},...,y_{r}\}\subset \mathcal{L}.$ Thus,
\begin{eqnarray*}
f(x) &=&t(Y)=F(\delta _{Y})=F\left( \sum_{i=1}^{r}\delta _{y_{i}}\right)
=F^{\prime }\left( \sum_{i=1}^{r}\delta _{y_{i}}\right) = \\
\sum_{i=1}^{r}F^{\prime }(\delta _{y_{i}})
&=&\sum_{i=1}^{r}g_{i}(y_{i})=\sum_{i=1}^{r}g_{i}(h_{i}(x)).
\end{eqnarray*}
\end{proof}

\textbf{Definition 2.3.} \textit{A closed path $p=\{x_{1},...,x_{n}\}$ is
said to be minimal if $p$ does not contain any closed path as its proper
subset.}

\bigskip

For example, the set
$l=\{(0,0,0),~(0,0,1),~(0,1,0),~(1,0,0),~(1,1,1)\}$ considered
above is a minimal closed path with respect to the functions $
h_{i}(z_{1},z_{2},z_{3})=z_{i},~i=1,2,3.$ Adding the point
$(0,1,1)$ to $l$, we will have a closed path, but not minimal. The
vector $\lambda $ associated with $l\cup \{(0,1,1)\}$ can be taken
as $(3,-1,-1,-2,2,-1).$

A minimal closed path $p=\{x_{1},...,x_{n}\}$ has the following obvious
properties:

\begin{description}
\item[(a)] \textit{The vector $\lambda $ associated with $p$ by Eq. (2.1) is
unique up to multiplication by a constant;}

\item[(b)] \textit{If in (2.1), $\sum_{j=1}^{n}\left\vert \lambda
_{j}\right\vert =1, $ then all the numbers $\lambda _{j},~j=1,...,n,$ are
rational.}
\end{description}

Thus, a minimal closed path $p$ uniquely (up to a sign) defines the
functional

\begin{equation*}
~G_{p}(f)=\sum_{j=1}^{n}\lambda _{j}f(x_{j}),\text{ \
}\sum_{j=1}^{n}\left \vert \lambda _{j}\right\vert =1.
\end{equation*}

\bigskip

\textbf{Lemma 2.4.} \textit{The functional $G_{p,\lambda }$ is a linear
combination of functionals $G_{p_{1}},...,G_{p_{k}},$ where $p_{1},...,p_{k}$
are minimal closed paths in $p.$}

\begin{proof} Let $\left\langle p,\lambda \right\rangle $ be a closed path-vector pair of $
X$, where $p=\{x_{1},...,x_{n}\}$ and $\lambda =(\lambda _{1},...,\lambda
_{n})$. Assume that $p$ is not minimal. Let $p_{1}=$ $\{y_{1},...,y_{s}\}$
be a minimal closed path in $p$ and

\begin{equation*}
G_{p_{1}}(f)=\sum_{j=1}^{s}\nu _{j}f(y_{j}),\text{ }\sum_{j=1}^{s}\left\vert
\nu _{j}\right\vert =1.
\end{equation*}

To prove the lemma, it is enough to show that $G_{p,\lambda }$ is a linear
combination of $G_{p_{1}}$ and some functional $G_{l,\theta }$, where $l$ is
a closed path in $X$ and a proper subset of $p$. Without loss of generality,
we may assume that $y_{1}=x_{1}.$ Put

\begin{equation*}
t_{1}=\frac{\lambda _{1}}{\nu _{1}}.
\end{equation*}
Then the functional $G_{p,\lambda }-t_{1}G_{p_{1}}$ has the form

\begin{equation*}
G_{p,\lambda }-t_{1}G_{p_{1}}=\sum_{j=1}^{k}\theta _{j}f(z_{j}),
\end{equation*}
where $z_{j}\in p$, $\theta _{j}\neq 0$, $j=1,...,k$. Clearly, the
set $ l=\{z_{1},...,z_{k}\}$ is a closed path with the associated
vector $\theta =(\theta _{1},...,\theta _{k})$. Thus, we obtain
that $G_{l,\theta }=$ $ G_{p,\lambda }-t_{1}G_{p_{1}}$. Note that
since $x_{1}\notin l$, the closed path $l$ is a proper subset of
$p$.
\end{proof}

\textbf{Theorem 2.5.} \textit{1) Let $X$ have closed paths. A
function $ f:X\rightarrow \mathbb{R}$ belongs to the space
$\mathcal{B} (h_{1},...,h_{r};X)$ if and only if $G_{p}(f)=0$ for
any minimal closed path $p\subset X$ with respect to the functions
$h_{1},...,h_{r}$.}

\textit{2) Let $X$ has no closed paths. Then $\mathcal{B}
(h_{1},...,h_{r};X)=T(X).$}

\begin{proof} 1) The necessity is clear. Let us prove the sufficiency. On
the strength of Lemma 2.4, it is enough to prove that if
$G_{p,\lambda }(f)=0$ for any closed path-vector pair
$\left\langle p,\lambda \right\rangle $ of $ X $, then $f\in
\mathcal{B}(X).$

Consider a system of intervals $\{(a_{i},b_{i})\subset \mathbb{R}
\}_{i=1}^{r} $ such that $(a_{i},b_{i})\cap
(a_{j},b_{j})=\varnothing $ for all the indices $i,j\in
\{1,...,r\}$, $~i\neq j.$ For $i=1,...,r$, let $\tau _{i}$ be
one-to-one mappings of $\mathbb{R}$\ onto $(a_{i},b_{i}).$
Introduce the following functions on $X$:

\begin{equation*}
h_{i}^{^{\prime }}(x)=\tau _{i}(h_{i}(x)),\text{ }i=1,...,r.
\end{equation*}

It is clear that any closed path with respect to the functions $
h_{1},...,h_{r}$ is also a closed path with respect to the
functions $ h_{1}^{^{\prime }},...,h_{r}^{^{\prime }}$, and vice
versa. Besides, $ h_{i}^{\prime }(X)\cap h_{j}^{\prime
}(X)=\varnothing ,$ for all $i,j\in \{1,...,r\},~i\neq j.$ Then by
Lemma 2.2,

\begin{equation*}
f(x)=g_{1}^{\prime }(h_{1}^{\prime }(x))+\cdots +g_{r}^{\prime
}(h_{r}^{\prime }(x)),
\end{equation*}
where $g_{1}^{\prime },...,g_{r}^{\prime }$ are univariate functions
depending on $f$. From the last equality we obtain that

\begin{equation*}
f(x)=g_{1}^{\prime }(\tau _{1}(h_{1}(x)))+\cdots +g_{r}^{\prime }(\tau
_{r}(h_{r}(x)))=g_{1}(h_{1}(x))+\cdots +g_{r}(h_{r}(x)).
\end{equation*}
That is, $f\in \mathcal{B}(X)$.

2) Let $f:X\rightarrow \mathbb{R}$ be an arbitrary function. First
suppose that $h_{i}(X)\cap h_{j}(X)=\varnothing ,$ for all $i,j\in
\{1,...,r\}$,$ ~i\neq j.$ In this case, the proof is similar to
and even simpler than that of Lemma 2.2. Indeed, the set of all
$(\ast )$-points of $Y$ consists of a single point, since
otherwise we would have a closed path with two points, which
contradicts the hypothesis of the 2-nd part of our theorem.
Further, the well definition of the functional $F$ becomes
obvious, since the left-hand side of (2.3) also contradicts the
nonexistence of closed paths. Thus, as in the proof of Lemma 2.2,
we can extend $F$ to the space $\mathcal{ S}^{\prime }$ and then
obtain the desired representation for the function $f$ . Since $f$
is arbitrary, $T(X)=\mathcal{B}(X).$

Using the techniques from the proof of the 1-st part of our
theorem, one can easily generalize the above argument to the case
when the functions $ h_{1},...,h_{r}$ have arbitrary ranges.
\end{proof}

\textbf{Theorem 2.6.} \textit{$\mathcal{B}(h_{1},...,h_{r};X)=T(X)$ if and
only if $X $ has no closed paths.}

\begin{proof} The sufficiency immediately follows from Theorem 2.5. To
prove the necessity, assume that $X$ has a closed path
$p=\{x_{1},...,x_{n}\}$. Let $ \lambda =(\lambda _{1},...,\lambda
_{n})$ be a vector associated with $p$ by Eq. (2.1). Consider a
function $f_{0}$ on $X$ with the property: $ f_{0}(x_{i})=1,$ for
indices $i$ such that $\lambda _{i}\,>0$ and $ f_{0}(x_{i})=-1,$
for indices $i$ such that $\lambda _{i}\,<0$. For this function,
$G_{p,\lambda }(f_{0})\neq 0$. Then by Theorem 2.5, $f_{0}\notin
\mathcal{B}(X)$. Hence $\mathcal{B}(X)\neq T(X)$. The
contradiction shows that $X$ does not admit closed paths.
\end{proof}

The condition whether $X$ have closed paths or not, depends both
on $X$ and the functions $h_{1},...,h_{r}$. In the following
sections, we see that if $ h_{1},...,h_{r}$ are "nice" functions
(smooth functions with the simple structure. For example, ridge
functions) and $X\subset \mathbb{R}^{d}$ is a "rich" set (for
example, the set with interior points), then $X$ has always closed
paths. Thus the representability by linear combinations of
univariate functions with the fixed "nice" multivariate functions
requires at least that $X$ should not possess interior points. The
picture is quite different when the functions $h_{1},...,h_{r}$
are not "nice". Even in the case when they are continuous, we will
see that many sets of $\mathbb{R}^{d}$ (the unite cube, any
compact subset of that, or even the whole space $\mathbb{R} ^{d}$
itself) may have no closed paths. If disregard the continuity,
there exists even one function $h$ such that every multivariate
function is representable as $g\circ h$ over any subset of
$\mathbb{R}^{d}$.

\bigskip

\begin{center}
\bigskip {\large \textbf{3. The analogue of the Kolmogorov superposition
theorem for multivariate functions}}
\end{center}

Let $X$ be a set and $h_{i}:X\rightarrow \mathbb{R},$ $i=1,...,r,$ be
arbitrarily fixed functions. Consider a class $A(X)$ of functions on $X$
with the property: for any minimal closed path $p\subset X$ with respect to
the functions $h_{1},...,h_{r}$ (if it exists), there is a function $f_{0}$
in $A(X)$ such that $G_{p}(f_{0})\neq 0.$ Such classes will be called
"permissible" function classes. Clearly, $C(X)$ and $B(X)$ are both
permissible function classes (in case of $C(X),$ $X$ is considered to be a
normal topological space).

\bigskip

\textbf{Theorem 3.1.} \textit{Let $A(X)$ be a permissible function
class. If $\mathcal{B}(h_{1},...,h_{r};X)=A(X)$, then $\mathcal{B}
(h_{1},...,h_{r};X)=T(X).$}

\bigskip

The proof is simple and based on the material of the previous
section. Assume for a moment that $X$ admit a closed path $p.$ The
functional $G_{p}$ annihilates all members of the set
$B(h_{1},...,h_{r};X).$ By the above definition of permissible
function classes, $A(X)\ $contains a function $ f_{0}$ such that
$G_{p}(f_{0})\neq 0.$ Therefore, $f_{0}\notin
B(h_{1},...,h_{r};X)$. We see that the equality
$B(h_{1},...,h_{r};X)=A(X)$ is impossible if $X$ has a closed
path. Thus $X$ has no closed paths. Then by Theorem 2.6,
$\mathcal{B}(h_{1},...,h_{r};X)=T(X).$

\bigskip

\textbf{Remark.} In the "if part" of Theorem 3.1, instead of
$\mathcal{B} (h_{1},...,h_{r};X)$ and $A(X)$ one can take
$\mathcal{B} _{c}(h_{1},...,h_{r};X)$ and $C(X)$ (or
$\mathcal{B}_{b}(h_{1},...,h_{r};X)$ and $B(X)$) respectively.

\bigskip

The main advantage of Theorem 3.1 is that we need not check directly if the
set $X$ has no closed paths, which in many cases may turn out to be very
tedious task. Using this theorem, we can extend free-of-charge\ the existing
superposition theorems for classes $B(X)$ or $C(X)$ (or some other
permissible function classes) to all functions defined on $X.$ For example,
this theorem allows us to obtain the analogue of the Kolmogorov
superposition theorem for all multivariate functions defined on the unit
cube.

\bigskip

\textbf{Corollary 3.2.} \textit{Let $d\geq 2$ and $s_{q}$,
$q=1,...,2d+1,$ be the Kolmogorov functions (1.2). Then each
function $f:\mathbb{I} ^{d}\rightarrow \mathbb{R}$ can be
represented by the formula (1.3), where $ g_{q}$ are univariate
functions depending on $f.$}

\bigskip

It should be remarked that Sternfeld [31], in particular, obtained
that the formula (1.3) is valid for functions $f\in
B(\mathbb{I}^{d})$ provided that $ g_{q}$ are bounded functions
depending on $f$ (see [15, chapter 1] for more detailed
information and interesting discussions).

Let $X$ be a compact metric space and $h_{i}\in C(X)$, $i=1,...,r.$ The
result of Sternfeld (see Introduction) and Theorem 3.1 give us the
implications
\begin{equation*}
\mathcal{B}_{c}(h_{1},...,h_{r};X)=C(X)\Rightarrow \mathcal{B}
_{b}(h_{1},...,h_{r};X)=B(X)\Rightarrow
\mathcal{B}(h_{1},...,h_{r};X)=T(X). \eqno(3.1)
\end{equation*}

The first implication is invertible when $r=2$ (see [31]). We want to show
that the second is not invertible even in the case $r=2.$ The following
interesting example is due to Khavinson [15, p.67].

Let $X\subset \mathbb{R}^{2}$ consists of a broken line whose sides are
parallel to the coordinate axis and whose vertices are

\begin{equation*}
(0;0),(1;0),(1;1),(1+\frac{1}{2^{2}};1),(1+\frac{1}{2^{2}};1+\frac{1}{2^{2}}
),(1+\frac{1}{2^{2}}+\frac{1}{3^{2}};1+\frac{1}{2^{2}}),...
\end{equation*}

We add to this line the limit point of the vertices $(\frac{\pi
^{2}}{6}, \frac{\pi ^{2}}{6})$. Let $r=2$ and $h_{1},h_{2}$ be the
coordinate functions. Then the set $X$ has no closed paths with
respect to $h_{1}$ and $ h_{2}.$ By Theorem 2.6, every function
$f$ on $X$ is of the form $ g_{1}(x_{1})+g_{2}(x_{2})$,
$(x_{1},x_{2})\in X$. Now construct a function $ f_{0}$ on $X$ as
follows. On the link joining $(0;0)$ to $(1;0)$ $
f_{0}(x_{1},x_{2})$ continuously increases from $0$ to $1$; on the
link from $(1;0)$ to $(1;1)$ it continuously decreases from $1$ to
$0$; on the link from $(1;1)$ to $(1+\frac{1}{2^{2}};1)$ it
increases from $0$ to $\frac{1}{2} $; on the link from
$(1+\frac{1}{2^{2}};1)$ to $(1+\frac{1}{2^{2}};1+\frac{1
}{2^{2}})$ it decreases from $\frac{1}{2}$ to $0$; on the next
link it increases from $0$ to $\frac{1}{3}$, etc. At the point
$(\frac{\pi ^{2}}{6}, \frac{\pi ^{2}}{6})$ set the value of
$f_{0}$ equal to $0.$ Obviously, $ f_{0} $ is a continuous
functions and by the above argument, $
f_{0}(x_{1},x_{2})=g_{1}(x_{1})+g_{2}(x_{2}).$ But $g_{1}$ and
$g_{2}$ cannot be chosen as continuous functions, since they get
unbounded as $x_{1}$ and $x_{2}$ tends to $\frac{\pi ^{2}}{6}$.
Thus, $\mathcal{B} (h_{1},h_{2};X)=T(X)$, but at the same time
$\mathcal{B}_{c}(h_{1},h_{2};X) \neq C(X)$ (or, equivalently,
$\mathcal{B}_{b}(h_{1},h_{2};X)\neq B(X)$).

We have seen that the unit cube in $\mathbb{R}^{d}$ has no closed
paths with respect to some $2d+1$ continuous functions (namely,
the Kolmogorov functions $s_{q}$ (1.2)). From the result of
Ostrand [25] (see Introduction) it follows that $d$-dimensional
compact sets $X$ also have no closed paths with respect to some
$2d+1$ continuous functions on $X$. One may ask if there exists a
finite family of functions $\{h_{i}:\mathbb{R}^{d}\rightarrow
\mathbb{R}\}_{i=1}^{n}$ such that any subset of $\mathbb{R}^{d}$
does not admit closed paths with respect to this family? The
answer is positive. This follows from the result of Demko [7]:
there exist $2d+1$ continuous functions $\varphi _{1},...,\varphi
_{2d+1}$ defined on $\mathbb{R}^{d}$ such that every bounded
continuous function on $\mathbb{R}^{d}$ is expressible in the form
$\sum_{i=1}^{2d+1}g\circ \varphi _{i}$ for some $ g\in
C(\mathbb{R})$. This theorem together with Theorem 2.6 yield that
every function on $\mathbb{R}^{d}$ is expressible in the form $
\sum_{i=1}^{2d+1}g_{i}\circ \varphi _{i}$ for some
$g_{i}:\mathbb{R} \rightarrow \mathbb{R},~i=1,...,2d+1$. We do not
yet know if $g_{i}$ here or in Corollary 3.2 can be replaced by a
single univariate function. We also don't know if the number
$2d+1$ can be reduced so that the whole space of $ \mathbb{R}^{d}$
(or any $d$-dimensional compact subset of that, or at least the
unit cube $\mathbb{I}^{d}$) has no closed paths with respect to
some continuous functions $\varphi _{1},...,\varphi
_{k}:\mathbb{R} ^{d}\rightarrow \mathbb{R}$, where $k<2d+1$. One
of the basic results of Sternfeld [32] says that the dimension of
a compact metric space $X$ equals $ d$ if and only if there exist
functions $\varphi _{1},...,\varphi _{2d+1}\in C(X)$ such that
$\mathcal{B}_{c}(\varphi _{1},...,\varphi _{2d+1};X)=C(X)$ and for
any family $\{\psi _{i}\}_{i=1}^{k}\subset C(X),$ $k<2d+1$, we
have $ \mathcal{B}_{c}(\psi _{1},...,\psi _{k};X)\neq C(X).$ In
particular, from this result it follows that the number of terms
in the Kolmogorov superposition theorem cannot be reduced. But
since the equalities $\mathcal{B }_{c}(X)=C(X)$ and
$\mathcal{B}(X)=T(X)$ are not equivalent, the above question on
the nonexistence of closed paths in $\mathbb{R}^{d}$ with respect
to less than $2d+1$ continuous functions is far from trivial.

If disregard the continuity, one can construct even one function
$\varphi : \mathbb{R}^{d}\rightarrow \mathbb{R}$ such that the
whole space $\mathbb{R} ^{d}$ will not possess closed paths with
respect to $\varphi $ and therefore, every function
$f:\mathbb{R}^{d}\rightarrow \mathbb{R}$ will admit the
representation $f=g\circ \varphi $ with some univariate $g$
depending on $f$. Our argument easily follows from Theorem 2.6 and
the result of Sprecher [30]: for any natural number $d$, $d\geq
2$, there exist functions $h_{p}:\mathbb{I}\rightarrow
\mathbb{R}$, $p=1,...,d,$ such that every function $f\in
C(\mathbb{I}^{d})$ can be represented in the form

\begin{equation*}
f(x_{1},...,x_{d})=g\left( \sum_{p=1}^{d}h_{p}(x_{p})\right) ,
\end{equation*}
where $g$ is a univariate (generally discontinuous) function
depending on $f$ .

\bigskip

\begin{center}
\bigskip {\large \textbf{4. Ridge functions}}
\end{center}

The set $\mathcal{B}(h_{1},...,h_{r};X)$, where $h_{i}$, $i=1,...,r,$ are
linear functionals over $\mathbb{R}^{d}$, or more precisely, the set
\begin{equation*}
\mathcal{R}(X)=\mathcal{R}\left(
\mathbf{a}^{1},...,\mathbf{a}^{r};X\right) =\left\{
\sum\limits_{i=1}^{r}g_{i}\left( \mathbf{a}^{i}\cdot \mathbf{x}
\right) :\mathbf{x}\in X\subset
\mathbb{R}^{d},~g_{i}:\mathbb{R}\rightarrow
\mathbb{R},i=1,...,r\right\} \eqno(4.1)
\end{equation*}
appears in many works (see, for example, [1,2,8,12,14,19,20,27]).
Here, $ \mathbf{a}^{i}$, $i=1,...,r,$ are fixed vectors
(directions) in $\mathbb{R} ^{d}\backslash \{\mathbf{0}\}$ and
$\mathbf{a}^{i}\cdot \mathbf{x}$ stands for the usual inner
product of $\mathbf{a}^{i}$ and $\mathbf{x}$. The functions
$g_{i}\left( \mathbf{a}^{i}\cdot \mathbf{x}\right) $ involved in
(4.1) are \textit{ridge} functions. Such functions are used in the
theory of PDE (where they are called \textit{plane waves}, see,
e.g., [13]), in statistics (see, e.g., [3,11]), in computerized
tomography (see, e.g., [14,20]), in neural networks (see, e.g.,
[28] and a great deal of references therein). In modern
approximation theory, ridge functions are widely used to
approximate complicated multivariate functions (see, e.g.,
[4,6,10,18,19,22,23,24,26,27,34,36]). In this section, we are
going to make some remarks on the representation of multivariate
functions by sums of ridge functions.

The problem of representation of multivariate functions by
functions in $ \mathcal{R}\left(
\mathbf{a}^{1},...,\mathbf{a}^{r};X\right) $ is not new. Braess
and Pinkus [1] considered the partial case of this problem:
characterize a set of points $\left( \mathbf{x}^{1},...,\mathbf{x}
^{k}\right) \subset \mathbb{R}^{d}$ such that for any data
$\{\alpha _{1},...,\alpha _{k}\}\subset \mathbb{R}$ there exists a
function $g\in \mathcal{R}\left(
\mathbf{a}^{1},...,\mathbf{a}^{r};\mathbb{R}^{d}\right) $
satisfying $g(\mathbf{x}^{i})=\alpha _{i},$ $i=1,...,k$. In
connection with this problem, they introduced the notion of the
\textit{NI}-property (non interpolation property) and
\textit{MNI}-property (minimal non interpolation property) of a
finite set of points as follows:

Given directions $\{\mathbf{a}^{j}\}_{j=1}^{r}\subset \mathbb{R}
^{d}\backslash \{\mathbf{0}\}$, we say that a set of points
$\{\mathbf{x} ^{i}\}_{i=1}^{k}\subset \mathbb{R}^{d}$ has the
\textit{NI}-property with respect to
$\{\mathbf{a}^{j}\}_{j=1}^{r}$, if there exists $\{\alpha
_{i}\}_{i=1}^{k}\subset \mathbb{R}$ such that we cannot find a
function $ g\in \mathcal{R}\left(
\mathbf{a}^{1},...,\mathbf{a}^{r};\mathbb{R} ^{d}\right) $
satisfying $g(\mathbf{x}^{i})=\alpha _{i},$ $i=1,...,k$. We say
that \ the set $\{\mathbf{x}^{i}\}_{i=1}^{k}\subset
\mathbb{R}^{d}$ has the \textit{MNI}-property with respect to
$\{\mathbf{a}^{j}\}_{j=1}^{r}$, if $\{\mathbf{x}^{i}\}_{i=1}^{k}$
but no proper subset thereof has the \textit{ NI}-property.

It follows from Theorem 2.6 that a set
$\{\mathbf{x}^{i}\}_{i=1}^{k}$ has the \textit{NI}-property if and
only if $\{\mathbf{x}^{i}\}_{i=1}^{k}$ contains a closed path with
respect to the functions $h_{i}=\mathbf{a} ^{i}\cdot \mathbf{x},$
$i=1,...,r$ (or, simply, to the directions $\mathbf{a} ^{i},$
$i=1,...,r$) and the \textit{MNI}-property if and only if the set
$\{ \mathbf{x}^{i}\}_{i=1}^{k}$ itself is a minimal closed path
with respect to the given directions. Taking into account this
argument and Definitions 2.1 and 2.3, we obtain that the set
$\{\mathbf{x}^{i}\}_{i=1}^{k}$ has the \textit{NI}-property if and
only if there is a vector $\mathbf{m} =(m_{1},...,m_{k})\in
\mathbb{Z}^{k}\backslash \{\mathbf{0}\}$ such that

\begin{equation*}
\sum_{j=1}^{k}m_{j}g(\mathbf{a}^{i}\cdot \mathbf{x}^{j})=0,
\end{equation*}
for $i=1,...,r$ and all functions $g:\mathbb{R\rightarrow R}$. This set has
the \textit{MNI}-property if and only if the vector $\mathbf{m}$ has the
additional properties: it is unique up to multiplication by a constant and
all its components are different from zero. This special corollary of
Theorem 2.6 was proved in [1].

Since ridge functions are nice functions of simple structure, representation
of every multivariate function by linear combinations of such functions may
not be possible over many sets in $\mathbb{R}^{d}$. The following remark
indicates the class of sets having interior points.

\bigskip

\textbf{Remark.} Let $X\subset \mathbb{R}^{d}$ have nonempty
interior. Then $ \mathcal{R}\left(
\mathbf{a}^{1},...,\mathbf{a}^{r};X\right) \neq T(X).$

\bigskip

Indeed, let $\mathbf{y}$ be a point in the interior of $X$.
Consider vectors $\mathbf{b}^{i}$, $i=1,...,r,$ with sufficiently
small coordinates such that $\mathbf{a}^{i}\cdot
\mathbf{b}^{i}=0$, $i=1,...,r$. Note that the vectors $
\mathbf{b}^{i}$, $i=1,...,r,$ can be chosen pairwise linearly
independent. With each vector $\mathbf{\varepsilon }=(\varepsilon
_{1},...,\varepsilon _{r})$, $\varepsilon _{i}\in \{0,1\}$,
$i=1,...,r,$ we associate the point

\begin{equation*}
\mathbf{x}_{\mathbf{\varepsilon
}}=\mathbf{y+}\sum_{i=1}^{r}\varepsilon _{i} \mathbf{b}^{i}.
\end{equation*}

Since the coordinates of $\mathbf{b}^{i}$ are sufficiently small,
we may assume that all the points $\mathbf{x}_{\mathbf{\varepsilon
}}$ are in the interior of $X$. We correspond each point
$\mathbf{x}_{\mathbf{\varepsilon } } $ to the number
$(-1)^{\left\vert \mathbf{\varepsilon }\right\vert }$, where
$\left\vert \mathbf{\varepsilon }\right\vert =\varepsilon
_{1}+\cdots +\varepsilon _{r}.$ One may easily verify that the
pair $\left\langle \{ \mathbf{x}_{\mathbf{\varepsilon
}}\},\{(-1)^{\left\vert \mathbf{\varepsilon } \right\vert
}\}\right\rangle $ is a closed path-vector pair of $X$. Therefore,
by Theorem 2.6, $\mathcal{R}\left( \mathbf{a}^{1},...,\mathbf{a}
^{r};X\right) \neq T(X).$

It should be noted that the above method of construction of the
set $\{ \mathbf{x}_{\mathbf{\varepsilon }}\}$ is due to Lin and
Pinkus [19].

Let us now give some examples of sets over which the representation by
linear combinations of ridge functions is possible.

\begin{description}
\item[(1)] Let $r=2$ and $X$ be the union of two parallel lines not
perpendicular to the given directions $\mathbf{a}^{1}$ and
$\mathbf{a}^{2}$. Then $X$ has no closed paths with respect to
$\{\mathbf{a}^{1},\mathbf{a} ^{2}\}$. Therefore, by Theorem 2.6,
$\mathcal{R}\left( \mathbf{a}^{1}{,} \mathbf{a}^{2};X\right)
=T(X).$

\item[(2)] Let $r=2,$ $\mathbf{a}^{1}=(1,1)$, $\mathbf{a}^{2}=(1,-1)$ and $X$
be the graph of the function $y=\arcsin (\sin x)$. Then $X$ has no closed
paths and hence $\mathcal{R}\left( \mathbf{a}^{1}{,}\mathbf{a}^{2};X\right)
=T(X).$

\item[(3)] Let now given $r$ directions $\{\mathbf{a}^{j}\}_{j=1}^{r}$ and $
r+1$ points $\{\mathbf{x}^{i}\}_{i=1}^{r+1}\subset \mathbb{R}^{d}$
such that
\begin{eqnarray*}
\mathbf{a}^{1}\cdot \mathbf{x}^{i} &=&\mathbf{a}^{1}\cdot
\mathbf{x}^{j}\neq \mathbf{a}^{1}\cdot \mathbf{x}^{2}\text{, \ for
}1\leq i,j\leq r+1\text{, } i,j\neq 2 \\ \mathbf{a}^{2}\cdot
\mathbf{x}^{i} &=&\mathbf{a}^{2}\cdot \mathbf{x}^{j}\neq
\mathbf{a}^{2}\cdot \mathbf{x}^{3}\text{, \ for }1\leq i,j\leq
r+1\text{, } i,j\neq 3 \\
&&\mathbf{......................................} \\
\mathbf{a}^{r}\cdot \mathbf{x}^{i} &=&\mathbf{a}^{r}\cdot
\mathbf{x}^{j}\neq \mathbf{a}^{r}\cdot \mathbf{x}^{r+1}\text{, \
for }1\leq i,j\leq r.
\end{eqnarray*}
The simplest data realizing these equations are the basis
directions in $ \mathbb{R}^{d}$ and the points $(0,0,...,0)$,
$(1,0,...,0)$, $(0,1,...,0)$ ,..., $(0,0,...,1)$. From the first
equation we obtain that $\mathbf{x}^{2}$ cannot be a point of any
closed path in $X=\{\mathbf{x}^{1},...,\mathbf{x} ^{r+1}\}$.
Sequentially, from the second, third, ..., $r$-th equations it
follows that the points
$\mathbf{x}^{3},\mathbf{x}^{4},...,\mathbf{x}^{r+1}$ also cannot
be points of closed paths in $X$ respectively. Thus the set $X$
does not contain closed paths at all. By Theorem 2.6,
$\mathcal{R}\left( \mathbf{a}^{1},...,\mathbf{a}^{r};X\right)
=T(X).$

\item[(4)] Let given directions $\{\mathbf{a}^{j}\}_{j=1}^{r}$ and a curve $
\gamma $ in $\mathbb{R}^{d}$ such that for any $c\in \mathbb{R}$,
$\gamma $ has at most one common point with at least one of the
hyperplanes $\mathbf{a} ^{j}\cdot \mathbf{x}=c$, $j=1,...,r.$ By
Definition 2.1, the curve $\gamma $ has no closed paths and hence
$\mathcal{R}\left( \mathbf{a}^{1},...,\mathbf{a }^{r};\gamma
\right) =T(\gamma ).$
\end{description}

At the end we want to draw the reader's attention to one more
problem concerning the set $\mathcal{R}\left(
\mathbf{a}^{1},...,\mathbf{a} ^{r};X\right) $. The problem is to
determine if a given function $ f:X\rightarrow \mathbb{R}$ belongs
to this set. One solution is proposed by Theorem 2.5: consider all
minimal closed paths $p$ of $X$ and check if $ G_{p}(f)=0$. This
problem was considered by some other authors too. For example, Lin
and Pinkus [19] characterized the set $\mathcal{R}_{c}\left(
\mathbf{a}^{1},...,\mathbf{a}^{r};\mathbb{R}^{d}\right) $ in terms
of the ideal of polynomials vanishing at all points $\lambda
\mathbf{a}^{i}\in \mathbb{R}^{d},$ $i=1,...,r,$ $\lambda \in
\mathbb{R}$. Two more characterizations of $\mathcal{R}_{c}\left(
\mathbf{a}^{1},...,\mathbf{a} ^{r};\mathbb{R}^{d}\right) $ may be
found in Diaconis and Shahshahani [8].

\smallskip

\textbf{Acknowledgements.} This research was done during my stay at the
Technion - Israel Institute of Technology. I am grateful to Allan Pinkus for
helpful discussions, his comments and pointing out the papers [31-33] by
Sternfeld and the book [15] by Khavinson. I am also grateful to Vitaly
Maiorov for numerous fruitful discussions and his ideas on the subject of
the paper. Besides, my special thanks go to the anonymous referee for many
helpful comments and the idea of the proof of Lemma 2.4.

\end{document}